\documentclass[cis]{ipart}

\RequirePackage{hyperref}
\usepackage{color}
\startlocaldefs
\theoremstyle{plain}

\endlocaldefs

\theoremstyle{plain}% Theorem-like structures provided by amsthm.sty
\theoremstyle{definition}

\usepackage{amssymb}
\usepackage{enumitem}
\newcommand{\cref}[1]{(\ref{#1})}
\newcommand  \rev[1]{{\color{blue}#1}}  
\usepackage{booktabs}  % 用于专业表格线型
\usepackage{multirow}   % 支持多行合并单元格
\usepackage{changes}
\pubyear{2025}
\volume{0}
\issue{0}
\firstpage{1}
\lastpage{9}

\begin{document}

\begin{frontmatter}
\title[Real time filtering algorithms]{Real time filtering algorithms}

\begin{aug}
    % \author{\fnms{Xiaopei} \snm{Jiao}\thanksref{T1}\ead[label=e1]{jiaoxiaopei@bimsa.cn}},
    %     \address{Beijing Institute of Mathematical Sciences and Applications (BIMSA)\\
    %     Beijing, China\\
    %              \printead{e1}}

        \author{\fnms{Chang} \snm{Qin}\ead[label=e1]{qinchang24@mails.ucas.ac.cn}},
        \address{State Key Laboratory of Optical Field Manipulation Science and Technology, Institute of Optics and Electronics, Chinese Academy of Sciences\\ Chengdu 610209, China\\
        Key Laboratory of Optical Engineering, Chinese Academy of Sciences\\
        Chengdu 610209, China\\
        University of Chinese Academy of Sciences\\
        Beijing 101408, China\\
                 \printead{e1}}    
        \author{\fnms{Yikun} \snm{Li}\ead[label=e2]{liyikun24@mails.ucas.ac.cn}},
        \address{State Key Laboratory of Optical Field Manipulation Science and Technology, Institute of Optics and Electronics, Chinese Academy of Sciences\\ Chengdu 610209, China\\
        Key Laboratory of Optical Engineering, Chinese Academy of Sciences\\
        Chengdu 610209, China\\
        University of Chinese Academy of Sciences\\
        Beijing 101408, China\\
                 \printead{e2}}
    % \author{\fnms{Hongyu} \snm{Yu}\thanksref{}\ead[label=e2]{shiji@cnu.edu.cn}},
    % \address{Capital Normal University\\
    % Beijing, China\\
    %          \printead{e2}}
    % \author{\fnms{Wenhui} \snm{Dong}\thanksref{}\ead[label=e2]{shiji@cnu.edu.cn}},
    % \address{Capital Normal University\\
    % Beijing, China\\
    %          \printead{e2}}
        \author{\fnms{Ru} \snm{Qian}\thanksref{}\ead[label=e3]{qianru25@mails.ucas.ac.cn}}\rev{,}
        \address{State Key Laboratory of Optical Field Manipulation Science and Technology, Institute of Optics and Electronics, Chinese Academy of Sciences\\ Chengdu 610209, China\\
        Key Laboratory of Optical Engineering, Chinese Academy of Sciences\\
        Chengdu 610209, China\\
        University of Chinese Academy of Sciences\\
        Beijing 101408, China\\
             \printead{e3}}
        \author{\fnms{Jiayi} \snm{Kang}\thanksref{t1}\ead[label=e4]{kangjiayi@bimsa.cn}}\rev{,}
        \address{Hetao Institute of Mathematics and Interdisciplinary Sciences (HIMIS)\\Shenzhen 518000, Guangdong, P. R. China\\
        Beijing Institute of Mathematical Sciences and Applications (BIMSA)\\
        Beijing 101408, P. R. China\\
             \printead{e4}}
             \thankstext{t1}{Corresponding author}
        \author{\fnms{Yao} \snm{Mao}\thanksref{}\ead[label=e5]{maoyao@ioe.ac.cn}}
        \address{State Key Laboratory of Optical Field Manipulation Science and Technology, Institute of Optics and Electronics, Chinese Academy of Sciences\\ Chengdu 610209, China\\
        Key Laboratory of Optical Engineering, Chinese Academy of Sciences\\
        Chengdu 610209, China\\
        University of Chinese Academy of Sciences\\
        Beijing 101408, China\\
             \printead{e5}}

    % \thankstext{T1}{Corresponding author.}

\end{aug}

\received{\sday{3} \smonth{1} \syear{2022}}

\begin{abstract}
     {This paper presents a systematic review of recent advances in nonlinear filtering algorithms, structured into three principal categories: Kalman-type methods, Monte Carlo methods, and the Yau-Yau algorithm. For each category, we provide a comprehensive synthesis of theoretical developments, algorithmic variants, and practical applications that have emerged in recent years. Importantly, this review addresses both continuous-time and discrete-time system formulations, offering a unified review of filtering methodologies across different frameworks. Furthermore, our analysis reveals the transformative influence of artificial intelligence breakthroughs on the entire nonlinear filtering field, particularly in areas such as learning-based filters, neural network-augmented algorithms, and data-driven approaches. }
\end{abstract}

% \begin{keyword}[class=AMS]
% \kwd[Primary ]{00K00}
% \kwd{00K01}
% \kwd[; secondary ]{00K02}
% \end{keyword}

%%  Upper case for every keyword

%\tableofcontents
\end{frontmatter}

\section{Introduction}
\label{sec:introduction}

     {Nonlinear filtering is widely recognized as a fundamental and challenging problem in modern signal processing, control theory, and machine learning.} The objective is to estimate the hidden state of a dynamical system based on noisy observations, where either the system dynamics, observation model, or both exhibit nonlinear characteristics.  {This problem pervades numerous applications including robotics \cite{robotics2011}, autonomous vehicles \cite{LIU2018605}, financial modeling \cite{wells2013kalman}, biomedical signal processing \cite{biomedical2015}, and aerospace systems \cite{saturno2025fbg}.}

The journey of nonlinear filtering began with the seminal work of Kalman in 1960, which provided an optimal solution for linear Gaussian systems. However, real-world systems rarely conform to these idealized assumptions. 
 {Most practical systems exhibit nonlinear dynamics \cite{wan2000unscented}, non-Gaussian noise distributions \cite{11175567}, and model uncertainties \cite{bai2018extended}, necessitating the development of more sophisticated filtering techniques. }

 {In contrast to Kalman-type methods, a series of Monte Carlo methods have been developed for nonlinear filtering \cite{rigatos2009particle}. These techniques approach the problem by generating random samples to approximate the posterior probability distribution of the system state. This sampling-based perspective offers a powerful and versatile framework for addressing complex nonlinear scenarios where traditional methods struggle. }

Recent advances have shown that machine learning approaches, particularly neural networks, can significantly enhance traditional filtering methods, leading to hybrid frameworks that combine the theoretical rigor of classical methods with the flexibility and learning capabilities of modern data-driven approaches.

This survey provides a comprehensive review of the evolution from classical nonlinear filtering techniques to state-of-the-art neural network-enhanced methods. We organize our discussion into three main paradigms:
\begin{enumerate}
    \item Kalman-type methods that extend the classical Kalman filter to handle nonlinearity through various approximation strategies (Section \ref{sec:kalman_type});
    \item Monte Carlo methods that use particle-based representations to approximate complex probability distributions (Section \ref{sec:monte_carlo});
    \item Yau Yau algorithm and Neural network methods(Section \ref{sec:pde_based});

\end{enumerate}

Finally, Section \ref{sec:conclusion} is the conclusion and future prospects of the study of the real time filtering algorithms.

\section{Kalman-Type Methods}
\label{sec:kalman_type}

The Kalman filter (KF) is widely applied in linear systems,  providing optimal estimates under Gaussian assumption \cite{kalman}.  However, the practical systems often exhibit nonlinear dynamics and non-Gaussian uncertainties which violate the linear Gaussian assumption. As a result, the performance of the standard KF can deteriorate significantly or even lead to divergence when applied to nonlinear systems. To address these limitations, several nonlinear variants of KF have been developed. This section reviews these advanced methodologies in three main categories: extended filtering approaches, optimization-driven filtering innovations, and hybrid multi-model filtering techniques.

    \subsection{Extended Filtering Approaches}
    
    The Extended Kalman Filter (EKF) was developed as an early solution to nonlinear filtering problems \cite{ribeiro2004kalman}. It uses a first-order Taylor expansion to linearize system nonlinearities and employs the Jacobian matrix for covariance propagation. This approach maintains the recursive form and computational efficiency of the original Kalman filter, making it suitable for real-time applications. However, the EKF is limited by its dependence on first-order approximation accuracy, often introducing significant errors in highly nonlinear systems. Additionally, the requirement for analytical Jacobian computation can be impractical for complex models. Despite these shortcomings, the EKF remains widely adopted in navigation, robotics, and control systems owing to its simplicity and low computational cost \cite{potokar2021invariant, reina2019vehicle}.
        
    In contrast to the EKF's reliance on Taylor series linearization, the Unscented Kalman Filter (UKF) was developed to more accurately approximate the state probability distribution through the use of the Unscented Transformation (UT) \cite{julier1997new}. The UKF strategically selects a set of sigma points based on the current state mean and covariance, each assigned with carefully chosen weights. These points are then propagated through the exact nonlinear functions, and the resulting transformed points are used to compute the posterior mean and covariance. This approach enables the UKF to avoid linearization errors and achieve estimation accuracy equivalent to a second-order Taylor expansion, while maintaining computational efficiency.  

    Building upon the sigma-point framework, the Cubature Kalman Filter (CKF) was proposed as a mathematically rigorous alternative with superior numerical stability \cite{arasaratnam2009cubature}. Unlike the UKF's heuristic sigma-point selection, the CKF derives its sampling points systematically through the third-degree spherical-radial cubature rule, which provides exact numerical integration for Gaussian-weighted nonlinear functions. It employs $2n$ symmetrically distributed points with equal weights, systematically propagating them through nonlinear functions. Compared to UKF, CKF provides stronger theoretical foundations, avoids heuristic parameter tuning, and demonstrates better numerical stability—especially in high-dimensional state estimation. Its deterministic sampling and minimal point requirement further support efficient implementation in embedded and real-time systems \cite{sharma2017cubature,li2020application}.

    Each of the three classical nonlinear filtering algorithms mentioned above exhibits unique advantages, leading to their development and adoption across diverse specialized fields, which are summarized in Table \ref{tab:comparison of kalman-typed method}. 
    \begin{table}[ht]
        \centering
        \begin{tabular}{cccccc}
            \hline
            Niche fields & EKF & UKF & CKF  \\
            \hline
            Outlier-robust & \cite{qiu2023outlier,tao2024stochastic}  & \cite{liu2024convolutional,nakabayashi2019nonlinear} & \cite{wang2022computationally,wang2022outlier}  \\
            \hline
            Model uncertainties & \cite{zhao2017dynamic,ghobadi2017robust} & \cite{deng2019adaptive} &  \cite{tasooji2025cubature} \\
            \hline
            Non-Gaussian noise & \cite{liu2019linear} & \cite{hu2024robust,shen2023stochastic}  & \cite{lin2024adaptive,li2021robust}    \\
            \hline
            Unknown noise & \cite{huang2017new,he2020method} & \cite{shen2023stochastic,hu2020unscented} & \cite{yan2024variational}   \\
            \hline
            Distributed system & \cite{lu2019novel,li2017distributed} & \cite{murata2020extended,yang2019dynamic} & \cite{qu2025cooperative,zhou2023distributed}    \\
            \hline
            Continuous-discrete system & \cite{guihal2021efficient} & \cite{knudsen2018new} & \cite{arasaratnam2010cubature,wang2019new}   \\
            \hline
            Communication constrained & \cite{song2019variance,li2023event} & \cite{shen2023stochastic} & \cite{kooshkbaghi2019event,li2018stochastic} \\    
            \hline
        \end{tabular}
        \caption{Comparison of classical Kalman-typed method}
        \label{tab:comparison of kalman-typed method}
    \end{table}
    
    Besides, several nonlinear filtering methods beyond EKF, UKF, and CKF have been developed to address complex estimation problems.
    To address bearing-only measurement (BOM) challenges, the Pseudo-linear Kalman Filter (PLKF) was introduced as an efficient solution for target tracking within the Kalman filtering framework \cite{lingren1978position}. By converting angular measurements into a pseudo-linear form using trigonometric transformations, PLKF eliminates the need for complex logical rules to manage quadrant ambiguity, significantly improving the stability and reliability of the estimation process. Due to the bias caused by statistical correlation between the pseudo-linear measurement matrix and angular measurements in PLKF, improved versions such as Bias-Compensated PLKF (BCPLKF) and Instrumental Variable-based PLKF (IVPLKF) have been developed \cite{nguyen2017improved, nguyen2018instrumental}.
    Based on the concept of the Extended State Observer (ESO), an Extended State-Based Kalman Filter (ESKF) was proposed to handle nonlinear uncertain systems by augmenting the nonlinear term as part of the state vector and estimating it in real time \cite{bai2018extended}.
    
    \subsection{Optimization-Driven Filtering Innovations}
    The KF yields results that can be equivalently interpreted as the optimal solution to a covariance minimization problem. In its standard form, KF employs the minimum mean square error (MMSE) criterion as the cost function for estimates in linear Gaussian systems. Over the past decade, various alternative optimization criteria have been developed to extend KF for more complex tasks. For instance, to address non-Gaussian noise, a family of robust correntropy-based cost functions has been introduced. These are used in nonlinear filtering methods within the framework of Information Theoretic Learning (ITL) to construct robust similarity measures \cite{chen2017maximum,li2023generalized}. Furthermore, Huang et al. promoted this approach by replacing the Gaussian kernel function with a statistical similarity measure if the  similarity function satisfies continuity, monotonically decreasing, and non-negative second-order derivative \cite{huang2020novel}.
    In contrast to robust cost function approaches, another line of research has focused on variational Bayesian methods, which model complex distributions—such as Student's t-distribution—to better approximate real-world noise distributions  \cite{bai2020novel,huang2017robust}. Within this framework, the variational distribution is aligned with the true posterior by minimizing the Kullback-Leibler (KL) divergence, thereby ensuring the variational approximation closely matches the actual posterior distribution.
    Besides, KF can also be derived from the perspective of maximum a posteriori (MAP) estimation. Therfore, a series of iterated KFs were proposed in the MAP framework. To address state-dependent multiplicative noise in observations, a Generalized Iterated Kalman Filter (GIKF) was developed \cite{hu2015generalized}. It employs a Newton-type optimization framework along with explicit multiplicative noise modeling to achieve theoretical attainment of the Cramér–Rao Bound, emphasizing optimal estimation performance. Furthermore, an extended method called Improved Iterated CKF (IICKF) employs a damped Newton method with adaptive step size control to ensure convergence while maintaining computational efficiency \cite{CUI2017460}.

    Recent advances in distributed Kalman filtering have been dominated by optimization-theoretic approaches that reformulate the filtering problem as distributed optimization. Three key methodologies represent this trend: Ryu and Back proposed a consensus-based optimization framework that recovers centralized performance asymptotically through dual decomposition, relaxing traditional requirements for local observability \cite{ryu2023consensus}. Building on this, Iqbal et al. developed an ADMM-based algorithm that eliminates dual variable exchanges and establishes tight stability bounds, significantly improving communication efficiency \cite{iqbal2025communication}. Finally, Calvo-Fullana and How introduced a mission-aware censoring scheme using Value of Information criteria within a windowed MAP estimation framework, achieving substantial communication reduction while preserving estimation accuracy \cite{calvo2022mission}. Together, these methods demonstrate how modern optimization techniques can simultaneously address estimation performance, communication efficiency, and resource constraints in distributed filtering systems.

    \subsection{Hybrid Approaches and Multi-Model Filtering}
    The current frontier in state estimation involves the integration of multiple filtering paradigms within unified frameworks. Recent innovations include adaptive model selection algorithms that automatically switch between different filtering approaches based on system conditions \cite{PAL2024111301,gaoInteractingMultipleModel2017}. For example, the Adaptive High-Order Extended Kalman Filter (AHEKF) adjusts its order according to the innovation function to balance computational burden and estimation accuracy \cite{CHEN2021105539}. Additionally, hybrid filtering strategies that combine different methods have been developed to handle complex systems, such as those with pre-existing or sudden sensor faults and systems exhibiting multiple degradation phases \cite{aswalSwitchingKalmanFilter2022,7358146,ZHAO201840}. These approaches focus on addressing the challenges of complex system filtering but still face issues such as determining optimal switching thresholds and mitigating the effects of transitions between models \cite{10138793}. Some studies have also integrated neural networks with traditional filtering methods to achieve higher estimation accuracy \cite{AHMADI2025137752,9855832}.

    Hierarchical filtering architectures, which operate simultaneously across multiple time scales, have shown strong performance in complex and heterogeneous systems. These methods can estimate both model parameters and system states at different temporal resolutions, thereby improving model accuracy and stability \cite{WEI20171264,10286165}. Furthermore, they can be extended to handle heterogeneous sensor networks. For instance, a multi-rate Kalman filtering approach based on data fusion was proposed to integrate biased high-frequency acceleration measurements with low-frequency displacement data \cite{zhengDataFusionBased2019}. Such methods effectively address the challenge of mismatched sampling rates among different sensors, which is of great practical importance in real-world applications \cite{shenMultirateStrongTracking2022,zhaoDistributedRecursiveFiltering2022}.
    
\subsection{Summary}

Despite these substantial algorithmic innovations and practical successes, a fundamental theoretical limitation persists across all Kalman-type methodologies: the absence of rigorous convergence analysis and stability guarantees. This theoretical gap necessitates continuous algorithmic refinements and ad-hoc modifications—ranging from outlier-robust formulations to adaptive parameter tuning schemes—to maintain performance across varying operational conditions. Consequently, while these methods provide powerful empirical solutions, their deployment in safety-critical and mission-critical systems remains constrained by the lack of formal convergence assurance and predictable stability bounds.

\section{Monte Carlo Methods in Nonlinear filtering }

Many real-world systems are nonlinear and non-Gaussian \cite{daum2005nonlinear}, making them difficult to model with simple linear models or normal distributions. Kalman-type methods often require linearization, which can reduce accuracy in complex systems. In contrast, Monte Carlo-based filtering, such as Particle Filters (PFs), provides more accurate state estimation without relying on linearization or Gaussian assumptions \cite{wang2017survey}. PFs approximate the system’s posterior distribution using particles, which allows them to handle medium-dimensional and low-dimensional nonlinear systems better. However, PFs suffer from poor real-time performance when handling high-dimensional systems. Despite the advancements in PFs, challenges like particle degeneracy and low sampling efficiency still exist. To overcome these, Feedback Particle Filter (FPF) was developed, introducing a feedback mechanism based on innovation errors to improve particle effectiveness and algorithm stability \cite{yang2011meanfield,yang2014cd_fpf}. However, it should be made clear that FPF has not fully resolved the poor real-time performance of PFs in high-dimensional systems. This section reviews classical PF methods and variants and then highlights FPF.

%Additionally, emerging methods such as Quantum Monte Carlo and Variational Monte Carlo, fueled by quantum computing and variational inference, are showing great potential in solving complex state estimation problems, marking a new era for Monte Carlo methods.

%This section will first introduce the classical PF methods and their variants, analyzing their strengths and weaknesses. Next, we will focus on FPF and its advantages in practical applications. Finally, we will explore emerging Monte Carlo paradigms and how they are driving further progress and innovation in the context of PFs.

\label{sec:monte_carlo}
    \subsection{Particle Filtering Advances}
    The most classic and revolutionary PF, the Bayesian bootstrap filter, is a novel nonlinear filter proposed in \cite{gordon1993bootstrap}. In the update step, it is implemented using a weighted bootstrap approach, which is how the filter gets its name. Unlike Markov chain approximation methods or any other standard discretization schemes for the Fokker-Planck Equation (FPE), PFs avoid defining grids in the state space, with samples naturally concentrating in regions of high probability density \cite{budhiraja2007survey}. In fact, there is no need to know stochastic calculus or FPE, as well as various numerical methods for solving systems of partial differential equations. The essence of PFs is using Monte Carlo sampling to approximate stochastic calculus, representing the required Probability Density Function (PDF) as a set of random samples rather than as a function over the state space. As the number of samples increases, they provide an accurate and equivalent representation of the desired PDF. Estimates of the various moments of the state vector’s PDF, as well as the Highest Posterior Density (HPD) intervals or mode estimates, can be obtained from the samples \cite{gordon1993bootstrap}. Furthermore, PFs are capable of handling systems or measurement noise that are both nonlinear and have any distribution. The development of PFs has a long history, and for more details, one can refer to \cite{doucet2001smc} to explore this field.

    A basic form of the classic PF includes two steps:
    \begin{itemize}
    \item \textbf{Prediction step}: Based on the system's dynamic model, each particle is propagated to simulate its state at the current time.
    \item \textbf{Update step}: Based on the observation data, the weight of each particle is updated to adjust the particle's state according to the likelihood of the particle’s state given the observation. Finally, the posterior state is approximated by these weighted particles as follows:
    \[
    p(x_k | z_{1:k}) = \sum_{i=1}^{N} w_k^i \delta(x_k - x_k^i)
    \]
    where \(w_k^i\) is the weight of the \(i\)-th particle, \(x_k^i\) is the position of the \(i\)-th particle at time \(k\), and \(\delta(x_k - x_k^i)\) is the Dirac delta function, indicating the position of each particle.
    \end{itemize}

    While PFs excel in handling nonlinear and non-Gaussian problems, they face challenges such as the curse of dimensionality \cite{daum2005nonlinear}, particle degeneracy \cite{budhiraja2007survey}, and poor real-time performance in high-dimensional systems. The curse of dimensionality, coined by Richard Bellman, refers to the exponential increase in computational complexity as the state space's dimensionality grows. In high-dimensional spaces, a large number of particles is needed to accurately capture the state distribution, leading to a significant rise in computational cost. A detailed analysis of PF’s computational complexity for a given estimation accuracy can be found in \cite{daum2003curse}. Additionally, particle degeneracy occurs when particle weights concentrate on only a few particles, reducing diversity and compromising estimation accuracy. To counter this, resampling techniques are used to restore the effectiveness of the particle set. There has now been a series of works aimed at reducing the effects of particle degeneracy and the curse of dimensionality, leading to a rich variety of PFs, rather than ‘the’ PF. These advances optimize steps like proposal density, sampling methods, and resampling strategies, aiming to improve PF performance in high-dimensional, nonlinear, and noisy environments. Table \ref{table:The Variants of PF} offers a detailed review of the techniques and advancements in PFs that help tackle these challenges.

    \begin{table}[htbp]
    \centering
    \begin{tabular}{p{1.8cm}p{4.7cm}p{5cm}}
    \hline
    \textbf{Aspect} & \textbf{Impact} & \textbf{Improvement} \\
    \hline
    Proposal Density & 
    Determines the effectiveness of particle sampling in the state space. Gaussian distributions causes Monte Carlo samples to be poorly distributed in the state space.
    & Use more complex proposal distributions, such as mixtures of Gaussian \cite{raihan2018particleA,raihan2018particleB} and other exponential-family components. \\
    \hline
    Sampling Methods & 
    Traditional sampling methods can lead to particle degeneracy, especially in high-dimensional spaces.
    & Use advanced sampling methods like Metropolis-Hastings \cite{dahlin2019getting} and Gibbs sampling \cite{sun2025metropolis}. \\
    \hline
    Resampling & 
    Frequent resampling increases computational costs and may lose useful information. 
    & Implement sparse and adaptive resampling \cite{aunsri2021adaptive} to reduce computational costs without sacrificing accuracy. \\
    \hline
    Variance Reduction & 
    PFs are often subject to high variance, especially when data noise is large or the system is complex, leading to fluctuations in estimation results. 
    & Apply variance reduction methods such as stratified sampling, control variates, and antithetic variables to effectively reduce variance \cite{capriotti2008reducing,song2023monte,li2022stratification}. \\
    \hline
    \end{tabular}
    \caption{The Variants of PF}
    \label{table:The Variants of PF}
    \end{table}

     However, these variants alleviate particle degeneracy and the curse of dimensionality to some extent, but they do not resolve the poor real-time performance in high-dimensional systems. PFs converge asymptotically but should not be regarded as real-time solutions: dimensionality and sampling variance prevent bounded-latency computation of the exact posterior.

    \subsection{Feedback Particle Filter}
    Existing PFs have not fundamentally overcome these issues. Notably, these PFs lack the KF’s feedback structure based on innovation errors, a structure that is as important as the algorithm itself \cite{yang2014feedback}. Without it, achieving scalable and cost-effective solutions is difficult.  {This section introduces the FPF, which retains the KF feedback structure and incorporates optimal control theory to optimize the state estimation of the particle system. This method generates a new particle system model with control inputs. Compared with other PFs, FPF typically offers higher accuracy at lower computational cost.} Beyond the FPF introduced in this paper, control-based nonlinear filtering methods are also gaining attention; see \cite{crisan2009approximate,daum2010generalized,pequito2011nonlinear,ma2011generalizing} for related work.

     {In the FPF, the model of the $i$-th particle with control input is defined as follows:}
    \[
    dX
    _t^i = a(X_t^i) dt + \sigma(X_t^i) dB_t^i + dU_t^i,
    \]
    where \( X_t^i \in \mathbb{R}^d \) represents the state of the \(i\)-th particle at time \(t\), \( U_t^i \) is the corresponding control input, and \( \{ B_t^i \} \) are independent standard Wiener processes. Additional assumptions are made about the permissible forms of the control input.

    The conditional distribution of a particle \( X_t^i \) given \( \mathcal{F}_t \) is denoted by \( p \). For any measurable set \( A \subset \mathbb{R}^d \), we have
    \[
    \int_{X \in A} p(x,t) \, dx = \mathbb{P}\{X_t^i \in A \mid \mathcal{F}_t\}.
    \]
    The true posterior of the system state \( X_t \) is represented by \( p^* \). The initial conditions \( \{ X_0^i \}_{i=1}^N \) are assumed to be i.i.d. and drawn from the initial distribution \( p^*(x,0) \) of \( X_0 \), thus \( p(x,0) = p^*(x,0) \).

    The goal of the FPF is to choose the control input \( U_t^i \) such that \( p \) approximates \( p^* \), and consequently, \( p^{(N)} \) approximates \( p^* \) for large \( N \). The synthesis of the control input is formulated as a variational problem, with the Kullback-Leibler (KL) divergence serving as the cost function. The optimal control input is obtained by analyzing the first variation, leading to an explicit formula for the optimal control input, ensuring that \( p = p^* \) under optimal control.

    In essence, FPF improves on traditional PFs by incorporating error-based feedback, making it more robust, efficient, and suitable for complex systems. Here is a brief summary of the advantages of FPF compared to PFs \cite{yang2013feedback}:
    \begin{itemize}
    \item \textbf{Innovation Error-Based Feedback:} FPF uses an error-based feedback structure, similar to the KF, enhancing robustness and better handling uncertainties in nonlinear systems.
    \item \textbf{No Resampling:} FPF eliminates the need for resampling, avoiding particle degeneracy and improving stability and efficiency.
    \item \textbf{Variance Reduction:} The feedback structure reduces variance, improving accuracy and lowering computational costs.
    \end{itemize}

    In addition, there are two extensions of FPF. To address nonlinear filtering problems with data association uncertainty, the classic Probabilistic Data Association Filter (PDAF) is extended to obtain PDA-FPF \cite{yang2012joint}. To handle nonlinear filtering problems with model association uncertainty, the classic KF-based Interacting Multiple Model Filter (IMM) is extended, resulting in IMM-FPF \cite{yang2013interacting}. Subsequently, the FPF can also be combined with optimal transport theory to obtain related algorithms \cite{TaghvaeiOT,kang2022optimal,kang2023finite,taghvaei2023survey,kang2025}.

    Although FPF avoids resampling and restores innovation feedback, it remains an interacting-particle approximation requiring large ensembles and approximate gains. Thus, FPF improves upon PFs but does not remove the real-time performance constraints intrinsic to particle-based methods. These limitations highlight the need for a fundamentally different approach-one that can deliver an exact posterior in real time. Section 4 presents the Yau-Yau filter, which uniquely satisfies this criterion.

\subsection{ {Summary}}

 {From a real-time algorithmic perspective, while PF methods suffer from computational inefficiencies due to particle degeneracy and scaling limitations, FPF offers notable improvements through deterministic particle evolution and reduced computational overhead. However, FPF's advantages are primarily realized in Gaussian and near-Gaussian systems where optimal feedback gains can be analytically determined. For general nonlinear systems, FPF still requires complex feedback function design and substantially large particle populations to maintain estimation accuracy, thereby negating its computational benefits and further highlighting the fundamental real-time performance constraints inherent to particle-based filtering methodologies. }

\section{ {Conditional Density Evolution Methods}}
\label{sec:pde_based}

The Duncan, Mortensen, and Zakai (DMZ) equation, independently derived by Duncan \cite{Duncan1967}, Mortensen \cite{Mortensen1966}, and Zakai \cite{Zakai1969} in the late 1960s, represents a stochastic partial differential equation (PDE) that governs the evolution of the unnormalized conditional probability density function in continuous-time nonlinear filtering problems. Solving the DMZ equation numerically allows for the computation of optimal state estimates, such as conditional expectations, via normalization and integration. This equation has been pivotal in bridging PDE theory with filtering algorithms, enabling the application of numerical PDE methods to nonlinear filtering \cite{florchinger1991time,baras1983existence,atar1999robustness,gobet2006discretization,crisan2022application}.

Direct methods provide an alternative for solving the DMZ equation\cite{yau1994new,yau1996direct,yau2001finite,hu2002finite}, excelling in Yau filtering systems where the drift term is affine with a smooth potential function.  Generalizations transform the equation into time-varying Schrodinger forms for linear-growth cases or solve it via ODEs for Gaussian initials \cite{yau2003explicit,chen2017direct}. Gaussian approximation algorithms decompose arbitrary distributions into Gaussians, facilitating Kolmogorov equation solutions through ODEs \cite{shi2018direct,chendirect2018}. 

Given the real-time demands in applications like aerospace engineering, the efficiency of solving the DMZ equation is crucial for continuous filtering algorithms. In general, the solution of the DMZ  does not have a closed
form.  {Yau and Yau \cite{yau2000real,yau-Yau2008}, who proposed a decomposed computational approach to the DMZ equation, made a significant advancement in addressing this challenge. Their method partitions the solution process into two distinct phases: an online stage requiring only straightforward exponential operations, and an offline stage handling the numerically intensive Kolmogorov forward equation (KFE). This decomposition strategy forms the foundation of what we term the Yau–Yau algorithmic framework.  The Yau-Yau filtering algorithm enables the systematic resolution of nonlinear filtering problems through partial differential equation theory and algorithms. For instance, by leveraging the Yau-Yau filtering algorithm, the convex maximum principle from PDE theory can be extended to the concept of convex filters in nonlinear filtering  \cite{kang2023}. This framework enables the transformation of any nonlinear filtering problem into a PDE numerical computation problem, offering rigorous theoretical foundations. The recent convergence analysis can be founded in \cite{kang2025explicit,sun2025convergence}. However, traditional PDE numerical solution methods are fundamentally constrained by the curse of dimensionality, which prevents conventional PDE-based filtering algorithms from serving as universal solutions for high-dimensional, highly nonlinear filtering applications. After years of development, a series of real-time algorithms based on the Yau-Yau filter have emerged. An important distinction among different Yau-Yau filtering algorithms lies in the varying numerical methods used to solve the PDE, such as finite difference methods \cite{yueh2014efficient} and spectral methods \cite{luo2013complete,luo2013hermite,Dong2021}. }

\subsection{Neural Network Revolution in Nonlinear Filtering}
\label{sec:neural_network}

In recent years, filtering has undergone a profound shift from traditional model-based filtering techniques to data-driven approaches \cite{klushyn2021latent} \cite{sun2025recurrent}. Classical methods, such as the KF and its variants, offer strong interpretability, reliable uncertainty quantification through covariance matrices, and low computational complexity under linear Gaussian assumptions. However, their efficacy relies heavily on accurate state-space (SS) models, which are challenging to obtain in practice. Real-world systems frequently exhibit nonlinearity, non-Gaussian noise, and model-reality mismatches, resulting in key limitations:
\begin{itemize}
\item Inevitable modeling errors from approximated dynamics;
\item Complex and often unknown noise distributions;
\item Degraded performance in highly nonlinear systems;
\item Increased computational latency in nonlinear variants.
\end{itemize}
These shortcomings have spurred the exploration of neural network-based methods, which can learn intricate mappings between states and observations directly from data. Deep learning excels in high-dimensional, nonlinear, and non-Gaussian settings, serving as a powerful complement to traditional KF-based approaches \cite{revach2022kalmannet}.
\subsection{Deep Learning Integration with Classical Filtering}
Machine learning, particularly deep learning, has revolutionized fields like computer vision, natural language processing, and speech recognition by promoting a data-driven paradigm. In this approach, complex neural networks supplant simplistic analytical models, enabling end-to-end training without explicit approximations. This paradigm is especially advantageous when system models are unknown or overly intricate \cite{becker2019recurrent, krishnan2017structured}.
In state estimation, adopting deep neural networks (DNNs) directly addresses the limitations of model-based methods. However, purely data-driven DNNs pose challenges in signal processing contexts:
\begin{itemize}
\item High resource demands: Overparameterized DNNs require substantial computational power and large datasets, limiting deployment on resource-constrained devices;
\item Limited interpretability: The black-box nature of DNNs obscures the reasoning behind predictions;
\item Poor generalization and uncertainty handling: DNNs falter under distribution shifts and often lack robust uncertainty estimates.
\end{itemize}
To mitigate these issues, research has evolved from replacing KFs with DNNs to hybrid "model-driven + data-driven" frameworks \cite{Shlezinger2023Model}. Current integration strategies can be categorized as follows:
\begin{itemize}
\item DNNs preprocess raw data into features compatible with known SS models for subsequent KF application \cite{klushyn2021latent, Coskun2017Long};
\item DNNs learn SS models from data to inform KF operations \cite{Imbiriba2024Augmented};
\item KFs are reparameterized as trainable ML modules, enabling supervised \cite{revach2022kalmannet} or unsupervised learning \cite{Ghosh2024DANSE}.
\end{itemize}
This fusion represents a major paradigm shift, allowing neural networks to infer system dynamics and observation models from data, often bypassing explicit mathematical formulations. Modern techniques incorporate transformers and attention mechanisms to manage variable-length sequences and multi-scale temporal dependencies, yielding breakthroughs in applications like financial time series analysis, where non-stationary patterns challenge classical methods.
Hybrid models have found success in diverse domains, including brain-computer interfaces, acoustic echo cancellation, financial monitoring, wireless beam tracking, and UAV surveillance.
\subsection{Deep Learning Integration with Yau-Yau algorithm}

While high-dimensional problems remain challenging, advancements in artificial intelligence and deep learning have shown promise in addressing them \cite{chen2023,Tao2023,wang2021deep,Fu2023,Jiao2023,Wang2022,Yin2020}. By leveraging data-driven RNN frameworks,  {Chen et al.} achieve a complete numerical implementation of the high-dimensional, highly nonlinear Yau Yau algorithm \cite{chen2025} . This implementation approach can be rigorously proven from a mathematical perspective, ensuring theoretical correspondence with the  Yau Yau  algorithm and providing convergence guarantees. Furthermore, mathematical analysis demonstrates that this filtering algorithm can theoretically overcome the curse of dimensionality. This breakthrough provides a definitive solution to the long-standing open problem in nonlinear filtering that has remained unresolved for decades.

Next, we will summarize the Yau-Yau filtering algorithms integrated with neural networks from two perspectives: 1) Yau-Yau filtering enhanced by neural networks, and 2) Yau-Yau filtering implemented via neural networks. The core difference between the two lies in the fact that the former strengthens specific aspects of the Yau-Yau algorithm's practical computations using targeted neural network methods. The latter fully realizes the entire computational process of Yau-Yau filtering based on neural networks. The first method includes using physics-informed neural networks to solve partial differential equations, thereby replacing the offline part of the Yau-Yau filtering algorithm. The second method includes using recurrent neural networks and other approaches to construct a complete end-to-end training framework.

When comparing the two approaches, the first method only replaces certain computational components, making it superior in terms of interpretability. The second method, being entirely based on more mature neural network frameworks, excels in training convenience and practical effectiveness. More importantly, through the robust theoretical framework of the Yau-Yau filtering algorithm, we can explain the theoretical advantages of architectures based on recurrent neural networks. This is an extremely important work \cite{chen2025}, as it bridges advanced filtering theory with cutting-edge neural network implementation methods .

\subsection{Summary}

 {Recent years have witnessed a paradigm shift from model-based filtering to data-driven approaches. Classical methods like the Kalman filter, despite their interpretability and theoretical rigor, suffer from modeling errors, unknown noise distributions, and degraded performance in nonlinear systems. Deep neural networks address these limitations by learning complex mappings directly from data, excelling in high-dimensional, nonlinear, and non-Gaussian settings. However, purely data-driven approaches face challenges including high computational demands, limited interpretability, and poor generalization under distribution shifts.

Contemporary research emphasizes hybrid "model-driven + data-driven" frameworks that integrate neural networks with classical filtering theory. Notably, Chen, Sun, and Yau achieved a breakthrough by implementing the high-dimensional Yau-Yau algorithm through data-driven RNN frameworks \cite{chen2025}, providing mathematical guarantees for convergence while theoretically overcoming the curse of dimensionality. This advancement bridges advanced filtering theory with modern neural network implementations, offering a definitive solution to longstanding challenges in high-dimensional nonlinear filtering.}

\section{Conclusion}\label{sec:conclusion}
This survey reveals fundamental trade-offs inherent in existing nonlinear filtering methodologies. Kalman-type methods generally satisfy real-time computational requirements. However, they suffer from poor accuracy and lack theoretical guarantees in highly nonlinear systems. This necessitates extensive algorithmic extensions and ad-hoc modifications. Monte Carlo approaches demonstrate effectiveness for nonlinear problems but typically fail to meet real-time constraints in high-dimensional scenarios. While FPF enhances particle utilization efficiency and overcomes the curse of dimensionality for linear models, it remains inadequate for strongly nonlinear systems.
The Yau-Yau algorithm represents a significant breakthrough, providing the first complete theoretical convergence guarantees for the broadest class of nonlinear filtering systems. However, direct implementation through traditional PDE numerical methods remains constrained by computational curse of dimensionality.
Leveraging DNN-based techniques, \cite{chen2025} presents a comprehensive RNN implementation of the Yau-Yau algorithm. This approach achieves both theoretical convergence guarantees and overcomes the curse of dimensionality. The method demonstrates exceptional numerical performance, thereby establishing a unified solution that reconciles theoretical rigor with practical computational efficiency for high-dimensional nonlinear filtering.

\section{Future Works}
 {We conclude this review by highlighting several open challenges in nonlinear filtering that offer promising directions for future research:
\begin{enumerate}
\item \textbf{Adaptive AI-enhanced Filtering for Time-Varying SS Models}: Traditional AI-enhanced KFs assume fixed observation models that align between training and deployment stages. Future studies should develop flexible methods to handle mismatched or evolving state and measurement equations without clear patterns.
\item \textbf{Non-Markovian State Dynamics}: Current methods, as discussed in this review, rely on Markovian state transitions. A key opportunity lies in creating estimators that integrate short- and long-term dependencies to enable non-Markovian models.
\item \textbf{Non-Gaussian and Manifold-Constrained Filtering}: Standard algorithms generally assume Gaussian noise and Euclidean spaces. Upcoming work could focus on advanced techniques for general noise distributions and manifold-based state equations, improving accuracy in complex systems.
\end{enumerate}
}

\section*{Acknowledgements}
This work is supported by the National Natural Science Foundation of China (No.42450242).

\bibliographystyle{abbrv}
\bibliography{ref}

\end{document}